\documentclass[]{article}
\usepackage{amssymb}
\def\ba{\mathbf{a}}   
\def\bb{\mathbf{b}}   
\def\be{\mathbf{e}}   
\def\bi{\mathbf{i}}   
\def\bj{\mathbf{j}}   
\def\bk{\mathbf{k}}   
\def\bm{\mathbf{m}}   
\def\bx{\mathbf{x}}   
\def\by{\mathbf{y}}   

\def\C{\mathbb{C}} 
\def\G{\mathbb{G}} 
\def\H{\mathbb{H}}  

\def\R{\mathbb{R}}  

\def\no{\noindent}

\def\beq{\begin{equation}}
\def\eeq{\end{equation}}

\def\w{\wedge}

\usepackage{helvet}
\usepackage{courier}
\usepackage{graphicx}   
\usepackage{makeidx}
\usepackage{mathptmx}
\def\bpm{\begin{pmatrix}}
\def\epm{\end{pmatrix}}
\makeindex            
\begin{document}
\title{Spinors in Spacetime Algebra and Euclidean 4-Space}
\author{Garret Sobczyk
\\ Universidad de las Am\'ericas-Puebla
 \\ Departamento de F\'isico-Matem\'aticas
\\72820 Puebla, Pue., M\'exico
\\ http://www.garretstar.com}
\maketitle
\begin{abstract} This article explores the geometric algebra of Minkowski spacetime, and its relationship to the geometric algebra of Euclidean 4-space. Both of these geometric algebras are algebraically isomorphic to the $2 \times 2$ matrix algebra over Hamilton's famous quaternions, and provide the rich geometric framework for various important
topics in mathematics and physics, including stereographic projection and spinors, and both spherical and hyperbolic geometry.
In addition, by identifying the timelike Minkowski unit vector with the extra $4^{\rm th}$-dimension of Euclidean $4$-space, David Hestenes' Space-Time Algebra of Minkowski spacetime is unified with William Baylis' Algebra of Physical Space.

\smallskip
\no {\em AMS Subject Classification:} 15A66, 81P16
\smallskip

\no {\em Keywords:} geometric algebra, hyperbolic geometry,  
 Poincare disk, quaternions, Riemann sphere, spacetime algebra, spinors, stereographic projection . 

\end{abstract}

\section*{0\ \ \  \  Introduction}
 
 This work represents the culmination of the author's attempt to understand the geometric nature of $2$-component Pauli and $4$-component Dirac 
 spinors in terms of a fully equivalent geometrical theory. In \cite{S1/2}, our geometric theory is shown to be fully equivalent to the classical Dirac theory of a spin-$\frac{1}{2}$ particle, whereas the Dirac-Hestenes equation is only equivalent to the classical Dirac equation if the issue of parity is not taken into consideration \cite{H66}, \cite{pal10}, \cite{strub05}. In \cite{S2014}, the spinor algebra of $\C^2$ is formulated in terms of idempotents and nilpotents in the geometric algebra $\G_3$ of Euclidean $3$-space, including its representation on the Riemann sphere, and a new proof of the Heisenberg uncertainty principle. In \cite{S2015}, the ideas of the previous work are generalized to include the Dirac theory by introducing the concept of a $2$-component Pauli-like matrix over a 4-dimensional Abelian module spanned by the imaginary unit $i=\sqrt{-1}$ and $J:=i \gamma_{0123}$, sometimes referred to as $\gamma_5$.
 In \cite{Shopf2015}, \cite{Sgroup2016}, a general view is taken discussing a range of topics such as the notion of a geometric spinor, some basic ideas of relativity, the Hopf fibration, and the concept of a group manifold.
 
The viewpoint adopted in all of the above mentioned work is that the geometric concept of number makes possible many new insights that are otherwise lost in a tangle of different obsolete formalisms. The {\it geometrization} of the real number system is accomplished by the following 

\begin{quote}{\bf Axiom:}   The real number system can be geometrically extended to include new,
anti-commutative square roots of $\pm 1$, each new such square root representing
the direction of a unit vector along orthogonal coordinate axes of a Euclidean or
pseudo-Euclidean space.
\end{quote}

\no This basic concept suggests, for example, the notation $\G_3 := \R(\be_1,\be_2 ,\be_3)$ to define the geometric algebra $\G_3$ by extending the real number system $\R$ to include $3$ new anti-commuting square roots of $+1$, denoted by $\be_1,\be_2,\be_3$. 

Section 1, introduces the geometric algebra $\G_4$ of $4$-dimensional Euclidean space $\R^4$ in terms of $2\times 2$ matrices over the quaternions, the unit quaternions being bivectors in the subalgebra $\G_3$ of $\G_4$. Two {\it spectral bases} of $\G_4$  are introduced, and formulas for the transformation between them are derived. 

Section 2, introduces the fundamental formula for stereographic projection in $\R^4$ and $\R^{1,3}$, and relates the geometric algebras $\G_4$ of the $4$-dimensional Euclidean space $\R^4$ to the geometric algebra $\G_{1,3}$ of Minkowski spacetime $\R^{1,3}$. In \cite[p. 24]{H66}, Hestenes introduces the the {\it spitting} of the space-time of the Pauli algebra $\G_3$ into the even subalgebra of the geometric algebra $\G_{1,3}$ of Minkowski spacetime. By identifying the timelike vector $\gamma_0 \in \R^{1,3}$ with the vector $\be_0 \in \R^4$, essentially the same splitting is achieved in the geometric algebra $\G_4$. In this 
scenario, the spacelike Minkowski vectors $\gamma_k$ for $k=1,2,3$ become bivectors in $\G_4^2$. Except for a change in the metric, the fundamental equations governing stereographic projection in $\R^4$ and $\R^{1,3}$ are identical. When considered in
$\R^{1,3}$, stereographic projection projects the Poincare disk onto the positive unit hyperboloid of Minkowski spacetime, one of the
5 basic models of hyperbolic geometry \cite{canal97}.
 
 Section 3, studies classical $2$-component spinors defined in terms of their equivalent {\it geometric spinors} in the geometric algebras $\G_3$ and $\G_{1,2}$, sub-algebras of $\G_4$ and $\G_{1,3}$, respectively. The Dirac bra-ket formalism is introduced and formulas, first derived in \cite{S2014} for the probability of finding a spin $\frac{1}{2}$-particle in a given state, are derived on the
 Bloch sphere $S^2$, along with an analogous formula derived on the Minkowski hyperboloid $L^2$.
 
 Section 4, generalizes the idea of a $2$-component geometric spinor, defined in the previous section in the algebras
 $\G_3$ and $\G_{1,2}$, to the concept of a $2$-component {\it quaternion-valued} spinor in the geometric algebras $\G_4$ and $\G_{1,3}$. Quaternion-valued spinors have been considered by other authors, see for example \cite{leorod1998}, \cite{mori07}.
 
 Section 5, derives the relationship between the classical $4$-component Dirac spinors and the $2$-component quaternion spinors defined in Section 4. It also relates the $2$-component $E$-spinors, defined in \cite{S2015} directly to quaternion spinors.

\section{Quaternions and Euclidean 4-space}

The structure of the geometric algebra $\G_4$ of Euclidean $4$-space is beautifully expressed in terms
of $2\times 2$ matrices over the quaternions $\H$. Here, quaternions are identified with bivectors in $\G_3$,
\[ \bi = \be_{23}, \ \ \ \bj = \be_{13}, \ \ \  \bk = \be_{12}, \]
so that  
\[  \H := \R(\be_{23}, \be_{13}, \be_{12}) = \{q | \ \ q=x_0 + x_1 \be_{23}-x_2 \be_{13} + x_3 \be_{12}, \ {\rm for} \ 
x_0, x_1,x_2,x_3 \in \R \}. \]
The minus sign in the third term of $q$, allows us to write
 \beq q=x_0 + i \bx= x_0+I( \bx \w \be_0) \label{quaternionq} \eeq
  where $\bx = x_1 \be_1+x_2 \be_2+x_3 \be_3$ is the {\it position vector} in the geometric subalgebra
  $\G_3 :=\R(\be_1,\be_2, \be_3) \subset \G_4$,
with $i:=\be_{123}$ and $I:=\be_{0123}$. When working in $\G_4$, or in its sub-algebra $\G_3$, the operation of {\it reverse} the order of the product of vectors of an
element $g\in \G_4$ is always denoted by $g^\dagger$.

A {\it spectral basis} of $\G_4$ over the quaternions is
\beq  \pmatrix{1 \cr i} \be_+ \pmatrix{1 & -i} = \pmatrix{ \be_+ & - i \be_- \cr i \be_+  & \be_-}, \label{spectralbasise} \eeq
where the idempotents $\be_\pm := \frac{1}{2}(1 \pm \be_0)$ are defined in terms of the unit vector $\be_0\in \G_4$. 
 The matrices quaternion $[\be_\mu]$ for $\mu = 0,1,2,3$ are
\beq  [\be_0]= \pmatrix{1 & 0 \cr 0 & -1},  [\be_1]=\pmatrix{ 0 & \be_{23} \cr -\be_{23} & 0},  [\be_2] = \pmatrix{ 0 &- \be_{13} \cr \be_{13} & 0},  [\be_3]= \pmatrix{ 0 & \be_{12} \cr -\be_{12} & 0} . \label{quabasis} \eeq
This spectral basis for $\G_4$ is different than the basis introduced by P. Lounesto in \cite[p.86]{LP97}

The matrices of any of the other bases elements of $\G_4$ are then easily recovered by taking sums of products of the matrices $[\be_\mu]$. The matrix $[q]$ of the quaternion (\ref{quaternionq}) is particularly simple, as is the matrix $[i]$ of $i$. We have
\[   [q] = \pmatrix{q & 0 \cr 0 & q}, \quad {\rm and} \quad [i] = \pmatrix{0 & -1 \cr 1 & 0}. \]
We can now solve for the matrix $[\bx]$ of any position vector $\bx \in \R^3$,
\[  [\bx] = - [i]\big( [q]- [x_0]\big) = \pmatrix{ 0 & 1 \cr -1 & 0}\pmatrix{ i \bx & 0 \cr 0 &   i \bx} = \pmatrix{0 & i \bx \cr -i \bx & 0}. \]
A general vector $x:=x_0 \be_0+\bx \in \R^4$ has the matrix
\beq  [x]_e = \pmatrix{x_0 & i \bx \cr - i \bx & -x_0}, \label{matrixpositionvece} \eeq
with respect to the spectral basis (\ref{spectralbasise}) defined by the primitive idempotents $\be_\pm$.

Another, closely related quaternion spectral basis of $\G_4$ that was used in \cite{Sgroup2016}, is
\beq  \pmatrix{1 \cr \be_0} I_+ \pmatrix{1 & \be_0} = \pmatrix{ I_+ & \be_0 I_- \cr \be_0 I_+  & I_-}. \label{spectralbasisI} \eeq
The general vector $x=x_0 \be_0+\bx \in \R^4$ has the quaternion matrix
\beq  [x]_I = \pmatrix{0 & x_0+ i \bx \cr x_0 - i \bx & 0}, \label{matrixpositionvece} \eeq
with respect to this spectral basis defined by the idempotents $I_\pm$. It follows that 
   \[ x = \pmatrix{ 1 & i } \be_+ [x]_e \pmatrix{1 \cr -i} =  \pmatrix{ 1 & \be_0 } I_+ [x]_I \pmatrix{1 \cr \be_0 } . \]
   Multiplying this equation on the left and right by the appropriate column and row vectors, we find that
   \[ \pmatrix{1& i \cr- i & 1 }\be_+ [x]_e \pmatrix{1 & i \cr -i & 1 } = 
    \pmatrix{1& \be_0 \cr -i & \be_0 i }I_+ [x]_I \pmatrix{1& i \cr \be_0  & \be_0 i  } .  \]
 Next, we multiply the left and right sides of this result by $\be_+$, and simplify, to get
        \[ \be_+ [x]_e =\be_+ \pmatrix{1& 1 \cr -I & I }I_+ [x]_I \pmatrix{1 & -I  \cr 1  & I } \be_+ = 
  \be_+I_+ \be_+ \pmatrix{1& 1 \cr -1 & 1 } [x]_I  \pmatrix{1 & -1 \cr  1 & 1 } .  \]

   Taking the conjugate of this last equation by multipling on the left and right by $i$ and $-i$ then gives
       \[ \be_- [x]_e = 
  \be_- I_- \be_-  \pmatrix{1& 1 \cr -1 & 1 }[x]_I  \pmatrix{1 & -1 \cr  1 & 1 } .  \]
    Adding together these last two results, and noting that $\be_+ I_+ \be_++ \be_- I_- \be_- = \frac{1}{2}$, gives the desired relationship
    \beq   [x]_e = A [x]_I A^{-1}, \label{changetoI} \eeq
    where 
    \[  A= \frac{1}{\sqrt 2} \pmatrix{1 & 1  \cr- 1 & 1 } \quad {\rm and} \quad  A^{-1} = \frac{1}{\sqrt 2}  \pmatrix{1 & -1 \cr   1 & 1 } \]
  are unitary matrices.
    It follows that the matrix $[g]_I$ is related to the matrix $[g]_e$, for any element $g\in \G_4$, by the same equation.
    
     An interesting
    relationship between the spectral bases (\ref{spectralbasise}) and (\ref{spectralbasisI}) is given by
    \[   \pmatrix{ I_+ & \be_0 I_- \cr \be_0 I_+  & I_-} = \pmatrix{ i_-  \cr i_+} \pmatrix{\be_+ + i \be_- & \be_+ -i \be_-}. \]
             Even more simply, 
    \[  I_+ = 2 i_- \be_+ i_+ \quad \iff \quad \be_+ = 2 i_+ I_+ i_-,\]
  where	 $i_\pm : = \frac{1}{2}(1 \pm i)$  which can be used to derive the relationship
    \beq   \pmatrix{\be_+ & -i \be_- \cr i \be_+ & \be_-} = 2 \pmatrix{i_+ \cr - i_-} I_+ \pmatrix{i_- & - i_+}  = B
    \pmatrix{I_+ & \be_0 I_- \cr \be_0 I_+ & I_-} B^*  ,  \label{spec-change-basis} \eeq
    for the {\it singular} matrix
    \[ B = \frac{\sqrt 2}{2} \pmatrix{i_+ & i_- \cr - i_- & i_+} ,\]
    where $B^* = {\overline B}^T$. The change of the spectral basis formula (\ref{spec-change-basis}), involving the singular matrix $B$, is quite different than the change of the coordinates formula (\ref{changetoI}) with the invertible unitary matrix $A$. 
   
\section{Stereographic projection}

The geometric algebra $\G_4$ of Euclidean 4-space $\R^4$, and the
spacetime geometric algebra $\G_{1,3}$ of the {\it pseudo-Euclidean} space $\R^{1,3}$ are algebraically isomorphic. We have  
\beq  \G_4= \R(\be_0,\be_1,\be_2,\be_3) = \R(\gamma_0, \gamma_{10},
\gamma_{20}, \gamma_{30} ) =\R(\gamma_0, \gamma_1, \gamma_2, \gamma_3)=\G_{1,3} ,\label{iso1}  \eeq
where $ \be_0 \equiv \gamma_0$ {\rm and} $\be_k \equiv \gamma_{k0}$ 
for $k=1,2,3$. Except for a renaming of the generating elements of $\G_4$, they are the same as the elements in the 
geometric algebra $\G_{1,3}$ of Minkowski spacetime, also called {\it Space-Time Algebra} (STA), \cite{H66}. When working in the geometric algebra $\G_{1,3}$, or in its
sub-algebra $\G_{1,2}$, the operation of {\it reverse} of an element $g\in \G_{1,3}$ is always denoted by
$\tilde g$. Of course, this is a different operation than the operation of reverse in $\G_4$ denoted by $g^\dagger$ . 

Using the algebraic isomorphism (\ref{iso1}) to solve for the spacetime Dirac vectors in $\G_4$, we find that
\[ \be_0 = \gamma_0, \ \ \ {\rm and} \ \ \ \be_{k0} = \gamma_{k 0 } \gamma_0 = \gamma_k ,  \]
so that
\beq  \G_{1,3}=\R(\gamma_0, \gamma_1, \gamma_2, \gamma_3)= \R(\be_0,\be_{10},\be_{20},\be_{30}) =\R(\be_0,\be_1,\be_2,\be_3)=\G_{4}. \label{gammag4} \eeq
Note the peculiar role of the vector $\be_0 \in \G_4$, which is identified with the timelike Minkowski
vector $\gamma_0 \in \G_{1,3}$. Contrast this with the identification of the spacetime Minkowski
vectors $\gamma_k $ with the Euclidean bivectors $\be_{k0}$ for $k=1,2,3$.

All calculations in this paper are carried out in the geometric algebras $\G_4:=\G(\R^4)$ of the Euclidean
space $\R^4$, or in $\G_{1,3}:= \G(\R^{1,3})$ of the pseudo-Euclidean space $\R^{1,3}$. Since we have identified
the elements of $\G_{1,3}$ with elements in $\G_4$ in (\ref{iso1}) and (\ref{gammag4}), all calculations could be carried out exclusively in either of these geometric algebras, but for reasons of clarity we do not attempt this here. 

The basic equation governing  stereographic projection of the unit 3-sphere $S^3$ in $\R^4$ into the flat $xyz$-space of $\R^3$, is
\beq  \bm = \frac{2}{\hat \ba + \be_0},  \label{basicstereo} \eeq
where $\hat \ba \in S^3$. 
 Equation (\ref{basicstereo}) implies that $\bm \cdot \be_0 = 1$, since rationalizing the denominator gives
\[  \bm = \frac{2(\hat \ba + \be_0)}{(\hat \ba + \be_0)^2}= \frac{\hat \ba + \be_0}{1+ \be_0 \cdot \hat \ba},  \]
 so that
\beq \bm = \frac{2}{\hat \ba + \be_0} = \bx_m + \be_0 \label{basicstereo2}  \eeq
where $\bx_m\in \R^3$ is the stereographic projection from the point $-\be_0$. The point $\bx_m$ is also the 
orthogonal projection of $\bm$ along the $\be_0$ axis.

\begin{figure}
\begin{center}
\no\includegraphics[scale=.25]{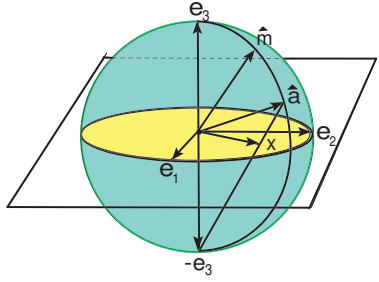}
\caption{Stereographic Projection from the South Pole at $-\be_3$ to the $xy$-plane. }
\label{sterox}
\end{center} 
  \end{figure}

The basic equation (\ref{basicstereo2}) further implies that
\beq   \hat \ba =\frac{\bm(2 - \bm \be_0)}{\bm^2}= \hat \bm \be_0 \hat \bm = e^{\theta \hat\bx_m \be_0} \be_0  , \label{defrotor} \eeq
where $ e^{\theta\hat\bx_m \be_0}:= (\hat \bm \be_0)^2$. The two sided {\it half-angle} form of this relationship is
\beq    \hat \ba =\hat \bm \be_0 \hat \bm =(\hat \bm \be_0) \be_0 (\be_0 \hat \bm) =  e^{\frac{\theta\hat\bx_m \be_0}{2}}\be_0
 e^{-\frac{\theta\hat\bx_m \be_0}{2}} .\label{defrotorhalf} \eeq 
The {\it rotor} $e^{\theta\hat\bx_m \be_0}$ rotates the North Pole at $\be_0 \in S^3$, $\theta$ radians into the point $\hat \ba \in S^3$. 
In Figure \ref{sterox}, we have replaced $\be_0$ with $\be_3$ and taken a cross-section of the Riemann $2$-sphere $S^2$ in $\R^3$,
in the plane of the bivector $\hat \bx_m \w \be_3$ through the origin. We see that the stereographic projection 
from the South pole at the point $-\be_3$, to the point $\hat \ba$ on the Riemann sphere,
passes through the point $\bx_m={\rm proj}(\bm)$ of the point $\bm$ onto the $xy$-plane through the origin with the normal vector $\be_3$.
 Stereographic projection is
 just one example of conformal mappings, which have important generalizations in higher dimensions \cite{Sob2012,SNF}. 

Returning to the geometric algebra $\G_4$, we
 expand $ e^{\theta\hat\bx_m \be_0}$ getting
 \beq   e^{\theta\hat\bx_m \be_0} = \cos \theta +\hat\bx_m \be_0 \sin \theta =  \frac{1-\bx_m^2+2\bx_m \be_0}{1+\bx_m^2},
  \label{expxm} \eeq
 so that
       \[            \cos \theta= \frac{1-\bx_m^2}{1+\bx_m^2}, \quad \sin \theta = \frac{2|\bx_m|}{1+\bx_m^2}, \]
     where $0\le \theta < \pi$ is the positive angle measured from the North Pole $\be_0$ to the point $\hat \ba \in S^3$.

     Using (\ref{defrotor}) and (\ref{expxm}), $\hat \ba = \frac{(1- \bx_m^2)\be_0 + 2 \bx_m}{1+\bx_m^2}$ which we use to calculate 
    \beq  d \hat \ba = \frac{2(1+ \bx_m^2)d \bx_m - 4 (\bx_m + \be_0) \bx_m \cdot d \bx_m}{(1+\bx_m^2)^2}
       \ \ {\rm and} \ \  (d \hat \ba)^2 =\frac{4 d\bx_m^2}{(1+\bx_m^2)^2},  \label{S-difarcmetric} \eeq
       giving the differential of arc and the Euclidean metric on $S^3$ in terms of the parameter $\bx_m \in \R^3$.  
 
 We have been exploring the geometry of the $3$-dimensional sphere $S^3$ in $\R^4$. Turning attention to the $3$-hyperboloid $L^3\subset \R^{1,3}$, defined by
\[   L^3= \{\hat a| \ \hat a^2 = 1 \} \]
where $\hat a =a_0 \gamma_0 + a_1 \gamma_1+a_2 \gamma_2 + a_3 \gamma_3$ for $a_\mu \in \R$, 
$\mu = 0,1,2,3$, and $a_0 \ge 1$.
Because of (\ref{gammag4}), it is possible to carry out all calculations in the geometric $\G_4$, but it is simpler to
work directly in $\G_{1,3}$.
  
By rewriting the  basic equation (\ref{basicstereo2}) in the form
\beq  m = \frac{2}{\hat a + \gamma_0} =x_m+\gamma_0,  \label{basicstereo3} \eeq
where now $x_m \in \R^{0,3}$, the same calculations for stereographic projection in $\R^4$ can be carried over to
 $\R^{1,3}$. All of the arguments used to derive (\ref{defrotor}),  (\ref{defrotorhalf}) and (\ref{expxm}) still apply, 
except in (\ref{expxm}) we need to change the trigonometric functions to hyperbolic trigonometric functions.
Thus,
\beq  m= \frac{2}{\hat a + \gamma_0}  = x_m+\gamma_0  \quad \iff \quad \hat a = \hat m \gamma_0 \hat m = e^{ \phi \hat \bx_m} \gamma_0 , \label{MA} \eeq
where $\hat \bx = \hat x_m \gamma_0$ is a vector in $\G_3^1$, or a bivector in $\G_{1,3}^2$, and
\beq e^{ \phi \hat \bx_m} = \cosh \phi +\hat \bx_m \sinh \phi = (\hat m \gamma_0)^2. \label{expphix} \eeq 
The two sided {\it half-angle} form of this relationship is
\[    \hat a =\hat m \gamma_0 \hat m =(\hat m \gamma_0) \gamma_0 (\gamma_0 \hat m) =  e^{\frac{\phi \hat \bx_m}{2}} \gamma_0
 e^{-\frac{\phi \hat \bx_m}{2}} .\] 
  \begin{figure}
\begin{center}
\no\includegraphics[scale=.35]{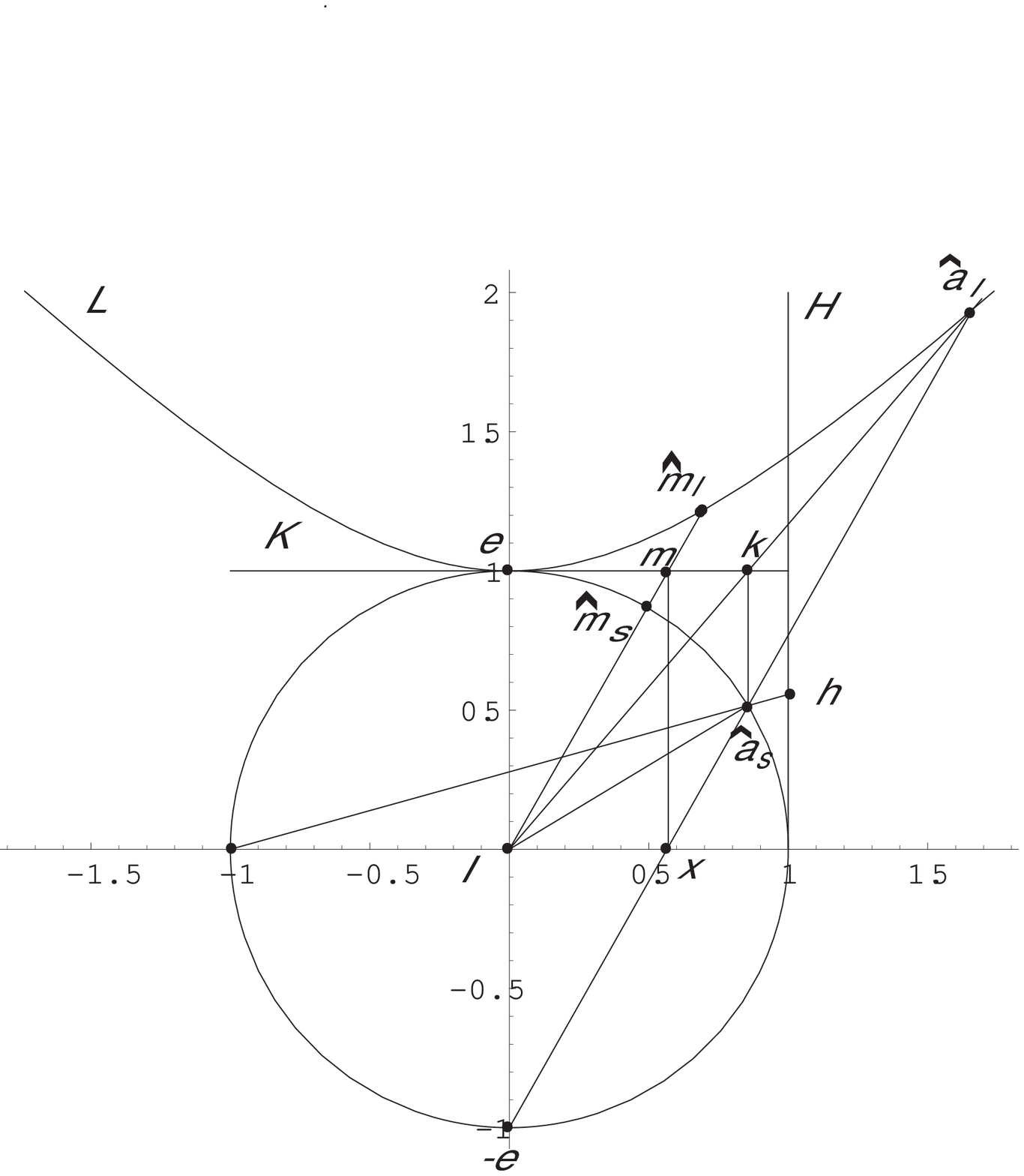}
\caption{Stereographic Projection on the ray from the South Pole at $-\be_0$ to the 3-plane at $x \in\R^{0,3}$, and onto the the hyperboloid at the point $\hat a$. We are making the identification $\be=\be_0=\gamma_0$.  Also shown are various well-known models of hyperbolic geometry, labeled by L, K, H, I and S. The idea for this Figure is borrowed from \cite{canal97}. }
\label{steroxL}
\end{center}
\end{figure}

\begin{figure}
\begin{center}
\no\includegraphics[scale=.60]{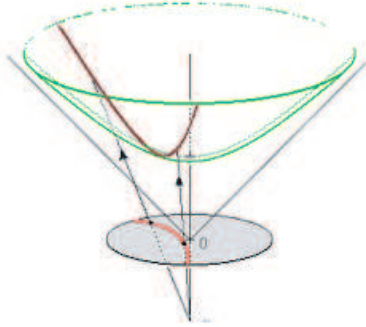}
\caption{Stereographic Projection from $-\gamma_0$, of the circular arc in the Poincare disk in $\R^{0,2}$, to the geodesic hyperbola on the hyperboloid $L^2$ in $\R^{1,2}$. (From Wikipedia: ``Poincare disk model''.) }
\label{steroxH}
\end{center}
\end{figure}

 Expanding $ e^{\phi \hat \bx_m}=(\hat m \gamma_0)^2$ in (\ref{expphix}), we find that
 \beq   e^{\phi \hat \bx_m} = \cosh \phi + \hat \bx_m \sinh \phi = (\hat m \gamma_0)^2 =
\frac{(1+\bx_m)^2 }{1-\bx_m^2}=\frac{1+\bx_m^2+2 \bx_m  }{1- \bx_m^2}, \label{expic} \eeq
     where $0\le \phi < \infty$ is the hyperbolic angle between $\gamma_0$ and $\hat a$  in $L^3$. 
It follows that
       \[            \cosh \phi= \frac{1+\bx_m^2}{1-\bx_m^2}, \quad  {\rm and} \quad \sinh \phi = \frac{2|\bx_m|}{1-\bx_m^2}.  \]
The rotor $e^{\phi \hat \bx_m}$ {\it boosts} $\gamma_0 \in L^3$, $\phi$ {\it hyperbolic} radians into $\hat a \in L^3$. Clearly, $\hat m$ is defined for all $x_m \in \R^{0,3}$ in the {\it open} unit $3$-sphere $x_m^2 < 1$. Figure \ref{steroxL} shows a cross-section of the $3$-hyperboloid,
taken in the plane of the bivector $\hat\bx=\hat x_m  \gamma_0$, through the origin. We see that the stereographic projection from the pole at the point $-\gamma_0$, to the point $\hat a$ on the hyperboloid,
passes through the point $x_m =proj(m)$ of the point $m$ onto the $xyz$-plane through the origin with the normal vector $\gamma_0$. Figure \ref{steroxH} shows the stereographic projection of the arc of a circle orthogonal to the circle bounding the Poincare disk into the corresponding geodesic hyperbola in $L^2$. 

   Using (\ref{MA}) and (\ref{expic}), $\hat a = \frac{(1+ \bx_m^2)\gamma_0 + 2 \bx_m \gamma_0}{1-\bx_m^2}$ which we use to calculate 
    \beq  d \hat a = \frac{2(1- \bx_m^2)d x_m + 4 (x_m + \gamma_0) \bx_m \cdot d \bx_m}{(1-\bx_m^2)^2},
       \  \  (d \hat a)^2 =\frac{-4 d \bx_m^2}{(1-\bx_m^2)^2}  \label{H-difarcmetric} \eeq
  giving the differential of arc and the hyperbolic metric on $L^3$ in terms of the parameter $\bx_m \in \R^{3}$, where $x_m:= \bx_m \gamma_0\in \R^{0,3}$ and $dx_m := d\bx_m \gamma_0\in \R^{0,3}$.  Note that the minus sign arises in (\ref{H-difarcmetric}) because of the hyperbolic metric on hyperboloid $L^3$. Contrast this with the closely related formula (\ref{S-difarcmetric}) for the sphere $S^3$.

\section{Geometric spinors}

We now restrict ourselves to geometric sub-algebras of $\G_4$ and $\G_{1,3}$. Consider
\[    \G_3 = \R(\be_1,\be_2,\be_3) \subset \G_4, \quad {\rm and} \quad \G_{1,2} = \R(\gamma_0,\gamma_1, \gamma_2)\subset \G_{1,3}. \]
The geometric algebra $\G_3$ is isomorphic to the Pauli algebra of matrices \cite{paulimat}, and is sometimes referred to as
the {\it Algebra of Physical Space} (APS), \cite{baysob2004}. In 1981, I called $\G_3$ ``Spacetime Vector Analysis'', because of its
close relationship to Gibbs-Heaviside vector analysis, \cite{S81}.
By a {\it geometric spinor}, or {\it g-spinor}, in the geometric sub-algebras $\G_3$ or $\G_{1,2}$, we mean the quantities
\[   \alpha := (\alpha_0+\alpha_1 \be_1)u_+\in \G_3 , \quad {\rm or} \quad    \alpha := (\alpha_0+\alpha_1 \gamma_1 )v_+\in \G_{1,2},   \]
respectively, 
where  $u_\pm: =\frac{1}{2}(1\pm\be_3)$ and $v_\pm = \frac{1}{2}(1 \pm \gamma_0)$, and $\alpha_0, \alpha_1$ are
{\it complex numbers} in $\G_3^{0+3}$ or $\G_{1,2}^{0+3}$, respectively.
A g-spinor $\alpha$ in $\G_3$ or $\G_{1,2}$ corresponds to a $2$-component spinor and a singular $(2 \times 2)$-matrix,  
\[ \alpha \ \longleftrightarrow \pmatrix{\alpha_0 \cr \alpha_1} \ \longleftrightarrow \  
 [\alpha]=\pmatrix{\alpha_0 & 0 \cr \alpha_1 & 0 }  . \]
 
For the $g$-spinor $\alpha \in \G_3$, we have 
 \beq   \alpha := \pmatrix{1 & \be_1 }u_+ [\alpha]\pmatrix{1 \cr \be_1}, \label{spinorg3} \eeq
 as explained in detail in \cite[(30)]{S2014}.  
 In the case that $\alpha \in \G_{1,2}$, the equation (\ref{spinorg3})
is modified to
 \beq   \alpha := \pmatrix{1 & \gamma_1 }v_+ [\alpha] \pmatrix{1 \cr -\gamma_1},  \label{spinorg12} \eeq
 taking into account that  $\gamma_1 \in\G_{1,2}$ has square $-1$.
The corresponding {\it bra-ket g-spinors} are defined by
\beq  |\alpha \rangle := \sqrt 2 \alpha, \quad {\rm and} \quad  \langle \alpha| :=\sqrt 2 \widetilde\alpha, \label{bra-ketdef} \eeq
 the factor $\sqrt 2$ being introduced for convenience. Whether a g-spinor $\alpha$ is in $\G_3$, or in
 $\G_{1,2}$, is made clear by context.

Factoring out $\alpha_0$, for a g-spinor $\alpha \in\G_3$, we find that
\beq \alpha = \alpha_0 (1+ \frac{\alpha_1}{\alpha_0}\be_1)u_+=
\rho e^{i \theta}\hat \bm u_+=\rho e^{i\theta}\hat \ba_+\hat \bm, \label{cangspinor3}  \eeq
where $i:=\be_{123}$, 
\[ e^{i\theta}:=\frac{\alpha_0}{\sqrt{\alpha_0 \alpha_0^\dagger}}, \quad {\rm and} \quad \rho:= \sqrt{\alpha_0 \alpha_0^\dagger + \alpha_1 \alpha_1^\dagger}=\sqrt{1+\bx_m^2}.\]
We also easily calculate the bra-ket relation
\beq   |\alpha \rangle \langle \alpha| = 2 \rho^2 \hat \bm u_+ \hat \bm = 2 \rho^2 \hat \ba_+ .  \label{hatarel} \eeq 

On the other hand, factoring out $\alpha_0$, for a g-spinor $\alpha \in\G_{1,2}$, we find that
\beq  \alpha = \alpha_0 (1+ \frac{\alpha_1}{\alpha_0}\gamma_1)v_+=
\rho e^{i \theta}\hat m v_+=\rho e^{i\theta}\hat a_+\hat m, \label{cangspinor12} \eeq
where $i:= \gamma_{012}$, 
\[ e^{i\theta}:=\frac{\alpha_0}{\sqrt{\alpha_0 \alpha_0^\dagger}}, \quad {\rm and} \quad \rho:= \sqrt{\alpha_0 \tilde \alpha_0 - \alpha_1 \tilde\alpha_1}= \sqrt{1+x_m^2},\]
and
\[  \hat m = \frac{x_m+\gamma_0}{\sqrt{1+x_m^2}}, \quad {\rm and} \quad \hat a_+ := \hat m v_+ \hat m, \]
for $x_m = x_1 \gamma_1+x_2 \gamma_2 $. For the bra-ket spinor $|\alpha \rangle$, the calculation (\ref{hatarel})
gives the similar result
\beq   |\alpha \rangle \langle \alpha| = 2 \rho^2 \hat m v_+ \hat m = 2 \rho^2 \hat a_+ .  \label{hataboost} \eeq

The $\pm$ sign changes that occur in the definition of $\rho$ in (\ref{cangspinor3}) and (\ref{cangspinor12}), are a result of the change of signature in going from $\G_3$ to $\G_{1,2}$. 
The inner products $\langle \alpha | \beta \rangle$ between the g-spinors $\alpha$ and $\beta$ in $\G_3$ or $\G_{1,2}$
are defined by
\beq  \langle \alpha | \beta \rangle :=  2 \langle \alpha^\dagger \beta \rangle_{0+3}=\big \langle \langle \alpha | | \beta \rangle \big \rangle_{0+3}, \ \ {\rm or} \ \  \langle \alpha | \beta \rangle :=  2 \langle \tilde \alpha \beta \rangle_{0+3}=\big \langle \langle \alpha | | \beta \rangle \big \rangle_{0+3}
 \label{innerproductforspinors} \eeq
and the {\it norm squared} 
   \[  |\alpha|^2 :=\langle \alpha | \alpha \rangle = \alpha_0^\dagger \alpha_0 +
       \alpha_1^\dagger \alpha_1>0  \ \ {\rm or} \ \  |\alpha|^2 :=\langle \alpha | \alpha \rangle = \tilde\alpha_0 \alpha_0 -
     \tilde  \alpha_1 \alpha_1>0,   \]
     respectively.
       
 For $\alpha=\rho_a e^{i \theta_a} \hat \bm_a u_+$ and 
$\beta=\rho_b e^{i \theta_b} \hat \bm_b u_+$ in $\G_3$,
\beq  \langle \alpha | \beta \rangle =2 \rho_a \rho_b e^{i(\theta_b-\theta_a)}\langle u_+\bm_a\bm_b u_+\rangle_{0+3} .\label{inprodspin} \eeq
 If a spin $\frac{1}{2}$-particle is prepared in 
a unit Pauli g-state $\alpha= \hat \bm_a u_+$, what is the probability of finding it in a unit Pauli g-state $\beta= \hat \bm_b u_+$ immediately
thereafter?
Following equations (13) and (14) in \cite{Shopf2015},
\[  \langle \beta |  \alpha \rangle \langle \alpha | \beta \rangle=4\Big\langle (\alpha^\dagger \beta)^\dagger (\alpha^\dagger \beta) \Big\rangle_{0+3} \]
\[=4\Big\langle u_+ \hat \bm_b \hat \bm_a u_+ \hat \bm_a \hat \bm_b u_+ \Big \rangle_{0+3}  \] 
\[ =4\langle \hat \bm_b \hat \bb_+ \hat \ba_+ \hat \bb_+ \hat \bm_b \rangle_{0+3} \]
\beq=\langle(1+\hat \ba \cdot \hat \bb)u_+\rangle_{0+3}=
\frac{1}{2}(1+\hat \ba \cdot \hat \bb). \label{probadotb} \eeq

This relationship can be directly expressed in terms of $\bm_a $ and $\bm_b$. 
For $\hat \ba = \frac{2}{\bm_a}-\be_0$ and $\hat \bb=\frac{2}{\bm_b}-\be_0$, a short calculation
gives the result
\beq  \langle \beta |  \alpha \rangle \langle \alpha | \beta \rangle= 1- \frac{(\bm_a-\bm_b)^2}{\bm_a^2 \bm_b^2} = \frac{1}{2}(1+\hat \ba \cdot \hat \bb)
  \label{probadotb2}. \eeq
\[ \iff \quad  \frac{(\bm_a-\bm_b)^2}{\bm_a^2 \bm_b^2} = \frac{1}{2}(1-\hat \ba \cdot \hat \bb).\]
Interpreting this result in $\G_3$ shows that the probability of finding the particle in that Pauli g-state $\beta$
is directly related to the Euclidean distance between the points $\bm_a$ and $\bm_b$.

 When $\hat \bb = -\hat \ba$,  
\[  \frac{(\bm_a-\bm_b)^2}{\bm_a^2 \bm_b^2} = 1. \]
This occurs when $\bm_b := \frac{1}{\bm_a \w \be_0}\bm_a$, in which case
\[  \hat \bb = \hat \bm_b \be_0 \hat \bm_b =
-\hat \bm_a \be_0 \hat \bm_a=-\hat \ba \ \ {\rm and} \ \  \bm_b \cdot \bm_a = 0,  \] 
and
\[ \bm_b=\bx_b+\be_0=\frac{1}{\bx_a \be_0}(\bx_a+\be_0) \]
\beq =\frac{\be_0 \bx_a}{\bx_a^2}(\bx_a+\be_0)=-\frac{1}{\bx_a}+\be_0 .\label{perpma} \eeq
This means that a spin $\frac{1}{2}$-particle prepared in the state $\bm_a$ will have a zero probability of being
found in the state $\bm_b$, for a measurement taken immediately afterwards. It is worthwhile mentioning that
these ideas can be naturally generalized to Dirac spinors \cite{S2015}, \cite{S1/2}.

 If, instead, we are doing the calculation in $\G_{1,2}$, we find that the calculations (\ref{probadotb}) and (\ref{probadotb2}) remain  
 valid, so that
  \beq  \langle \beta |  \alpha \rangle \langle \alpha | \beta \rangle= 1- \frac{(m_a-m_b)^2}{m_a^2 m_b^2} =\frac{1}{2}(1+\hat a \cdot \hat b) \ge 1, \label{hypadotb} \eeq
 which clearly cannot represent a probability. The fact that hyperbolic dot product $\hat a \cdot \hat b$ can be expressed in terms of the hyperbolic difference $(m_a-m_b)^2$ for the points $m_a$ and $m_b$ on the hyperboloid $L^2$, suggest the name
 {\it Bloch hyperboloid} in analogy for the case of the {\it Bloch sphere} represented in (\ref{probadotb2}). In analogy to the
 Hopf Fibration for the Bloch sphere, we can similarly define a Hopf Fibration for the case of the Bloch hyperboloid, \cite{Shopf2015}.
  
  \section{Quaternion spinors}

     The g-spinors $\alpha \in \G_3$ and $\alpha \in \G_{1,2}$, defined in (\ref{spinorg3}) and (\ref{spinorg12}), can be generalized to a quaternion-valued spinor in either $\G_4$ or in $\G_{1,3}$, depending upon in which algebra we are working.
   Recalling (\ref{quaternionq}) and (\ref{iso1}), a quaternion in either $\G_4$ or $\G_{1,3}$ can be written in the form $q =x_0+  i\bx$ for $i=\be_{123}=\gamma_{0123}$.
   Just as for a g-spinor, a {\it quaternion spinor} 
   \[  \alpha_q := (q_0 + q_1 i)v_+ \quad \longleftrightarrow \quad \pmatrix{q_0 \cr q_1 } \quad \longleftrightarrow \quad
    [\alpha_q]:=  \pmatrix{q_0 & 0 \cr q_1 & 0}     . \] 
With respect to the spectral basis,
      \[   \pmatrix{1 \cr i }v_+ \pmatrix{1 & -i} = \pmatrix{v_+ &- i v_- \cr i v_+ & v_-}, \]
    which was first defined in (\ref{spectralbasise}) for $v_\pm = \frac{1}{2}(1\pm \be)$, with $\be=\be_0 = \gamma_0$, 
     the quaternion spinor $\alpha_q$ satisfies the equation 
     \[  \alpha_q =  \pmatrix{1 &  i }v_+ [\alpha_q]  \pmatrix{1 \cr -i}   . \]

    If $q_0 = x_0 + i \bx$ and $q_1 = y_0 + i \by$, then
   \[  q_0^\dagger q_1  = (x_0 - i \bx)(y_0 + i \by) =( x_0 y_0+\bx \cdot \by) +i( x_0 \by - y_0  \bx +  \bx \times \by) =q_0^\dagger \circ q_1 + q_0^\dagger \otimes q_1 ,  \]
 where
   \[ q_0^\dagger \circ q_1 :=\frac{1}{2}(q_0^\dagger q_1 + q_1 q_0^\dagger) = x_0 y_0+\bx \cdot \by+(x_0 \by-y_0 \bx)i \] 
   and 
    \[ q_0^\dagger \otimes q_1:=\frac{1}{2}(q_0^\dagger q_1 - q_1 q_0^\dagger)  = ( \bx \times \by) i . \]
        Also note that
   \[ \langle q_1 q_0^\dagger \rangle_0 = \frac{1}{2} (q_1 q_0^\dagger + q_0 q_1^\dagger) = x_0 y_0 + \bx \cdot \by \]
and 
  \[    \langle q_1 q_0^\dagger \rangle_1 =\frac{1}{2} (q_1 q_0^\dagger - q_0 q_1^\dagger) =i(x_0 \by - y_0 \bx - \bx \times \by) . \]
       
  For $\alpha_q=(q_0+q_1i)v_+$, then
  \beq \alpha_q = \frac{1}{q_0 ^\dagger}\big( q_0^\dagger q_0 + q_0^\dagger q_1 i \big)v_+ 
  = q_0 \Big[ \be_0 +\frac{ \langle q_0^\dagger  q_1\rangle_0}{q_0^\dagger q_0}i+\frac{\langle q_0^\dagger  q_1\rangle_2 }{q_0^\dagger q_0} i \be_0   \Big] v_+ =   \sqrt{ q_0 q_0^\dagger}  e^{\theta i \hat \bx} M v_+,   \label{Mcanform} \eeq  
    where 
   \[  e^{\theta i \hat \bx}:= \frac{q_0}{\sqrt{q_0 q_0^\dagger}}, \ \ {\rm and} \ \
    M =\Big[ \be_0 +\frac{ \langle q_0^\dagger  q_1\rangle_0}{q_0^\dagger q_0}i
 +\frac{ \langle q_0^\dagger  q_1\rangle_2}{q_0^\dagger q_0} i \be_0   \Big] .\]
 Expressing $M$ as an element of $\G_{1,3}$ gives
 \[ M = \gamma_0+\frac{ y_0 x - x_0 y - \gamma_{123}x\w y}{x_0^2-x^2} +\frac{ \gamma_{0123}(x_0y_0 - x\cdot y)}{x_0^2 - x^2} ,\]
 where $\bx = x \gamma_0$, $\by = y \gamma_0$ for $x,y \in \G_{0,3}^1$, and $i=\gamma_{0123}$.

 Noting that
  \[ M^2=1 - \frac{q_0^\dagger q_1}{q_0^\dagger q_0 }\Big[\frac{ \langle q_0^\dagger  q_1\rangle_2}{q_0^\dagger q_0}  -
  \frac{\langle q_0^\dagger  q_1\rangle_0}{q_0^\dagger q_0}\Big]
  =  1+ \frac{\langle q_0^\dagger  q_1\rangle_2^2- \langle q_0^\dagger  q_1\rangle_0^2}{(q_0^\dagger q_0)^2}=1-\frac{q_1^\dagger q_1}{q_0^\dagger q_0} ,  \]
  we arrive at the {\it canonical form}
\beq    \alpha_q = \rho    e^{\theta i \hat \bx} \hat M v_+,   \label{gencanformqspinor} \eeq
 where 
 \[  \rho =\sqrt{q_0^\dagger q_0}|M| =
              \sqrt{q_0^\dagger q_0 -q_1^\dagger q_1}  \quad {\rm and} \quad \hat M = \frac{M}{|M|}\]
    for the quaternion spinor $\alpha_q$. We say that $\alpha_q$ is a {\it unit quaternion spinor} if $\rho = 1$.
          Naturally, generalizing (\ref{bra-ketdef}), we define the bra-ket quaternion spinors
    \beq  | \alpha_q \rangle := \sqrt 2 \alpha_q, \quad {\rm and} \quad \langle \alpha_q | = \sqrt 2 \alpha_q^\dagger. \label{bra-ketquaternion} \eeq    
    
   Defining, and carrying out the same calculations for quaternion spinors that we did in (\ref{innerproductforspinors}), (\ref{inprodspin}), and (\ref{probadotb}), gives      
  \[  \langle \alpha_q | \beta_q \rangle = 2 \rho_a \rho_b \langle v_+ \hat M_a e^{-\theta_a i \hat \bx_a} 
   e^{\theta_b i \hat \bx_b} \hat M_b v_+ \rangle_{0+3} = 2 \rho_a \rho_b \langle e^{-\theta_a i \hat \bx_a }v_+ \hat M_a^\prime 
 \hat M_b^\prime v_+ e^{\theta_b i \hat \bx_b }\rangle_{0+3} ,    \]
   where $\hat M_a^\prime =e^{\theta_a i \hat \bx_a}   \hat M_a    e^{-\theta_a i \hat \bx_a }$ and 
    $\hat M_b^\prime =e^{\theta_b i \hat \bx_b}   \hat M_b    e^{-\theta_b i \hat \bx_b }$, and
    \[    \langle \beta_q | \alpha_q \rangle    \langle \alpha_q | \beta_q \rangle  =
                      4 \langle (\alpha_q^\dagger \beta_q)^\dagger ( \alpha_q^\dagger \beta_q) \rangle_{0+3}
                        =  4 \langle v_+ \hat M_b^\prime \hat B^\prime_+ \hat A_+^\prime \hat M_b^\prime v_+ \rangle_{0+3}    \]
                    \[ =       4 \langle v_+ \hat M_b^\prime \hat B^\prime_+ \hat A_+^\prime \hat M_b^\prime v_+\rangle_{0+3}
                    =    4 \langle  \hat M_b^\prime \hat B^\prime_+ \hat A_+^\prime \hat B_+^\prime \hat M_b^\prime \rangle_{0+3} 
                    = \frac{1}{2}\big(1 + \hat A^\prime \circ \hat B^\prime\big),\]
        where $\hat A^\prime = \hat M_a^\prime \gamma_0 \hat M_a^\prime$ and             
  $\hat B^\prime = \hat M_b^\prime \gamma_0 \hat M_b^\prime$. 
 
     A quaternion spinor $\alpha_q$ is said to be {\it orthogonal} if $  \langle q_0^\dagger q_1 \rangle_0 = 0$.
    For the orthogonal quaternion spinor $\alpha_q$, the canonical forms (\ref{Mcanform}),  (\ref{gencanformqspinor}) simplify to
    \beq  \alpha_q   =  \sqrt{ q_0 q_0^\dagger} e^{\theta i \hat \bx} \big( 1+\bx_m \big)
    \be_0 v_+ =   \sqrt{ q_0 q_0^\dagger} e^{\theta i \hat \bx} \big(\gamma_0+x_m \big)
     v_+,\label{orthogqspinor} \eeq
    and 
     \beq  \alpha_q  = \rho e^{\theta i \hat \bx}\hat  M v_+ =\rho  e^{\theta i \hat \bx} 
      \frac{(1+\bx_m)\be_0}{\sqrt{1- \bx_m^2}}  v_+,\label{orthogqspinor2} \eeq
    respectively, where   
  \[  \bx_m :=  \frac{y_0 \bx - x_0 \by - \bx \times \by}{x_0^2 + \bx^2}= \frac{y_0 (  x\,  \gamma_0) - x_0 (y \, \gamma_0) -i (x \w y)}{x_0^2 + \bx^2}   . \]
 We also have
     \[ |M| = \sqrt{M^2} =\sqrt{ 1 - \bx_m^2} = \sqrt{ 1 - \frac{y_0^2 + \by^2}{x_0^2 + \bx^2}}. \]
     
    When $\alpha_q$ is an orthogonal quaternion g-spinor, then with the help of (\ref{orthogqspinor2}), 
    \[ |\alpha_q \rangle \langle \alpha_q | = \Big[1 + (x_0^2 -  y_0^2 + \bx^2 -  \by^2)\be_0 -2 (x_0 \by - y_0 \bx - \bx \times \by)\Big]. \]
   In agreement with (\ref{orthogqspinor}), the canonical form 
   (\ref{cangspinor12}) for a g-spinor $\alpha \in \G_{1,2}$ generalizes to
   \[ \alpha_q = q_0 M v_+ = \rho e^{\theta i \hat \bx  } \hat M v_+, \]
    for an orthogonal quaternion g-spinor $\alpha_q \in \G_{1,3}$, 
    where 
    \[\bx_m=i q_0^{-1} q_1, \ \  M=(1+\bx_m)\be_0, \ \  \rho := \sqrt{q_0 q_0^\dagger - q_1 q_1^\dagger}, \ \
    {\rm and} \ \ e^{i \hat \bx \theta}=\frac{q_0}{\sqrt{ q_0 q_0^\dagger}}.  \]
    In this case, $ \langle \alpha_q | |\alpha_q \rangle = 2 \rho^2 v_+$ for $\rho = \sqrt{q_0 q_0^\dagger - q_1 q_1^\dagger}$.
    
    \section{Classical Dirac spinors}
    
    In earlier work \cite[(5)]{S2015}, the classical $4$-component Dirac spinors were considered by introducing the
    spectral basis
  \beq \pmatrix{1 \cr \be_{13} \cr \be_3 \cr \be_1} u_{++}  \pmatrix{1 & -\be_{13} & \be_3 & \be_1}
= \pmatrix{u_{++} & -\be_{13}u_{+-} & \be_3 u_{-+} & \be_1 u_{--}
 \cr \be_{13} u_{++} &u_{+-} & \be_1 u_{-+} & -\be_3 u_{--} \cr
 \be_3 u_{++} & \be_{1}u_{+-} & u_{-+} & -\be_{13} u_{--} \cr 
  \be_1 u_{++} & - \be_{3}u_{+-} & \be_{13} u_{-+} &  u_{--} }, \label{specbasisD} \eeq  
    for the primitive idempotents
  \[u_{\pm \mp}:=\frac{1}{4}(1\pm \gamma_0)(1 \mp i\gamma_{12})=\frac{1}{4}(1\pm \be_0)(1 \mp J \be_3), \]
  where $J:= -j i$ for $j:= \sqrt{-1}$ and $i=\be_{123} = \gamma_{0123}$. A geometric Dirac spinor is then defined by the
  correspondence
  \[ |\alpha \rangle_4 :=\pmatrix{\varphi_1 \cr \varphi_2 \cr \varphi_3 \cr \varphi_4} \ \ \longleftrightarrow \ \  
  \pmatrix{\varphi_1 & 0 & 0 & 0  \cr \varphi_2 & 0 & 0 & 0   \cr \varphi_3 & 0 & 0 & 0   \cr \varphi_4 & 0 & 0 & 0  }
  \ \ \longleftrightarrow \ \  |\alpha \rangle_g, \]
  where 
  \[  |\alpha \rangle_g :=( \varphi_1 + \varphi_2 \be_{13} + \varphi_3 \be_3 + \varphi_4 \be_1)u_{++},  \]
  and $\varphi_k := x_k + j y_k$ for $k=1,2,3,4$ and $x_k,y_k\in \R$.  
  
  As explained in detail in \cite[(11)]{S2015}, noting that $j u_{++}=\gamma_{21}u_{++} = i \be_3 u_{++}$, the geometric spinor $ |\alpha \rangle_g$ takes the form
\[ |\alpha \rangle_g =\Big( (x_1+x_4 \be_1 + y_4 \be_2+x_3 \be_3)+i(y_3+y_2 \be_1 - x_2 \be_2 + y_1 \be_3) \Big)u_{++}. \]
The geometric spinor $| \alpha \rangle_g$ can be rewritten in terms of the quaternions $q_0$ and $q_1$,
\[  |\alpha\rangle_g \ = (q_0 + i q_1)u_{++} =\Big( (x_0 +i \bx)+i(y_0+ i \by)\Big)v_+ E_+  \]
\[  = \Big(  (x_0 -\by )+i(y_0+ \bx)\Big)v_+ E_+,  \]
where $E_\pm := \frac{1}{2}(1 \pm J \be_3)$. It follows that the classical Dirac spinor $ | \alpha \rangle_4$ can be expressed in terms of the quaternion spinor $| \alpha \rangle_q$, by defining the components
\[  \varphi_1 = x_0 + ix_3, \ \varphi_2 = -x_2 + i x_1, \ \varphi_3 = -y_3+ iy_0, \ \varphi_4 = - y_1 -y_2 i . \]
It follows that the whole classical theory of $4$-component Dirac theory can be expressed equivalently in terms of
geometric quaternion spinors. 

\section*{Acknowledgement} I thank Universidad de Las Americas-Puebla for many years of support.

\end{document}